\documentclass{svproc}

\usepackage{url}

\usepackage{hyperref}
\hypersetup{
    colorlinks=true,
    linkcolor=blue,
    filecolor=blue,      
    urlcolor=blue,
    citecolor=blue}

\usepackage{amsmath}
\usepackage{amssymb} 

\def\iid{\overset{\textnormal{iid}}{\sim}} 

\makeatletter  
\@ifundefined{@currsizeindex} 
  {\RequirePackage{relsize}\let\dolarger\relsize} 
  {\def\dolarger#1{\larger[#1]}} 
\newcommand*\@@bigtimes[2]{\vphantom{\prod} 
  \vcenter{\hbox{\dolarger{4}$\m@th#1\mkern-2mu\times\mkern-2mu$}}} 
\newcommand*\bigtimes{\mathop{\mathpalette\@@bigtimes\relax}\displaylimits} 
\makeatother

\def\N{\mathbb{N}}\def\R{\mathbb{R}}\def\1{\mathbbm{1}}

\def\Fcal{\mathcal{F}}\def\Wcal{\mathcal{W}}

\def\cube{[0,1]^d}

\begin{document}
\mainmatter              
\title{Bayesian inference with Besov-Laplace priors for spatially inhomogeneous binary classification surfaces}
\titlerunning{Bayesian binary classification with Besov-Laplace priors}  
%
\author{Matteo Giordano}
\authorrunning{Matteo Giordano} 
%
%
\institute{ESOMAS Department, University of Turin, Corso Unione Sovietica 2018/bis, 10134 Turin, Italy\\
\email{matteo.giordano@unito.it}
}

\maketitle              

%
%
%
%
%

\begin{abstract}
In this article, we study the binary classification problem with supervised data, in the case where the covariate-to-probability-of-success map is possibly spatially inhomogeneous. We devise nonparametric Bayesian procedures with Besov-Laplace priors, which are prior distributions on function spaces routinely used in imaging and inverse problems in view of their useful edge-preserving and sparsity-promoting properties. Building on a recent line of work in the literature, we investigate the theoretical asymptotic recovery properties of the associated posterior distributions, and show that suitably tuned Besov-Laplace priors lead to minimax-optimal posterior contraction rates as the sample size increases, under the frequentist assumption that the data have been generated by a spatially inhomogeneous ground truth belonging to a Besov space.
\end{abstract}

%
%
%
%
%

\section{Introducion}

The binary classification problem consists in predicting a binary response $Y$ based on the value of a covariate $X$. Taking, without loss of generality, $Y$ to be a $0$-$1$ outcome, the canonical approach to this problem entails modelling $Y$, conditionally given $X$, as a Bernoulli random variable with probability of success $f(X)$. Given a labelled training dataset $(Y^{(n)},X^{(n)}):=\{(Y_i,X_i)\}_{i=1}^n$, the central task is then to estimate the `covariate-to-probability-of-success map' $x\mapsto f(x)$. Based on this, new unlabelled inputs $X_{n+1}$ can be classified according to the estimate of $f(X_{n+1})$, e.g.~to be of class $1$ if this probability exceeds $1/2$.

	In the present article, we are interested in the case where the classification surface $f$ is possibly spatially inhomogeneous, namely when it is generally flat or smooth across the domain, but it may exhibit localised high variation or even discontinuities. This is particularly important for a variety of applications in which there may be unknown `critical values' of the covariates that determine rapid increases or decreases in the probability of success.

	We shall consider the nonparametric Bayesian approach to this problem, endowing $f$ with prior distributions of Besov-Laplace type. These are popular priors on function spaces that were first introduced by Lassas et al.~\cite{LSS09}, and have since then enjoyed a widespread application in imaging and inverse problems, where they empirically achieve excellent performances in the recovery of spatially inhomogeneous objects of `bounded variation type', giving rise to edge-preserving and sparse reconstructions at the level of the maximum-a-posteriori (MAP) estimator. See e.g.~\cite{LP01,BD06,VLS09,SE15,KLSS23}, where many more references can be found.

	Besov-Laplace priors are constructed by modelling the random coefficients in a wavelet basis expansion via independent (scaled) Laplace distributions. This offers a `discretisation-invariant' alternative \cite{LSS09} to the famous `total-variation' prior of Rudin et al.~\cite{ROF92}, and furnishes a Bayesian model for functions in the Besov scale $B^\alpha_{11}, \ \alpha\ge0$. These spaces measure the local variability in an $L^1$-sense and describe the regularity by the $\ell^1$-decay of the wavelet coefficients, and are well-know to provide a mathematical characterisation for spatially inhomogeneous functions. See Section \ref{Sec:Notation} for details, and \cite[Section 1]{DJ98} for an overview on the relationship with the space of bounded variation functions. In contrast, the popular Gaussian priors, which entail $L^2$-type penalties, are suited to model Sobolev-regular functions with milder variations, and have been shown to perform sub-optimally in more structured recovery problems \cite{ADH21,giordano2022inability,agapiou2024laplace}.  We refer the reader to Appendix E.2 in \cite{agapiou2024heavy} for an empirical comparison of the performance of Besov-Laplace and Gaussian priors over simulated data in the white noise model with spatially inhomogeneous regression functions. Section 4.2.3 of the Ph.D.~thesis \cite{Savva2024} contains further extensive simulations results, and a nice visualisation of the superior reconstruction quality achieved by Besov-Laplace priors in the presence of spatially inhomogeneous ground truths, again in the white noise model. Alongside the aforementioned advantages, Besov-Laplace priors also maintain a favourable log-concave structure that enables computational approaches \cite{BG15} and analytical study \cite{ADH21}.

	Building on the breakthroughs in the theory of the frequentist analysis of nonparametric Bayesian procedures from the last two decades \cite{GGvdV00,shen2001rates,GvdV07,vdVvZ08,GN11}, Agapiou et al.~\cite{ADH21} have recently started a line of research on the performance of posterior inference with Besov-Laplace priors for the recovery of spatially inhomogeneous functions, in the large sample size limit and under the frequentist assumption that the data have been generated by a fixed `ground truth' belonging to a Besov space. Minimax-optimal posterior contraction rates in the white noise model have been obtained by \cite{ADH21}, later extended to drift estimation for multidimensional diffusions, density estimation and nonlinear inverse problems by \cite{GR22,giordano23besov,agapiou2024laplace}, respectively. Further, adaptive posterior contraction rates with related hierarchical procedures have been obtained by \cite{giordano23besov,agapiou2024adaptive}. We also mention the earlier related results by \cite{CN14,R13,AGR13}, and the recent examination by \cite{dolera2024strong}.

	In this paper, we shall study these issues for the binary classification problem, which to our knowledge has not been considered in the context of Besov-Laplace priors and inference for spatially inhomogeneous functions. In our main result, Theorem \ref{Theo:Main}, we shall show that for true classification surfaces belonging to the $B^\alpha_{11}$-scale, suitably tuned Besov-Laplace priors attain minimax-optimal posterior contraction rates. Our proof (in Section \ref{Sec:Proof}) is based on the general posterior contraction rate theory for independent and identically distributed (i.i.d.) binary classification data, e.g.~\cite[Chapter 2.5]{GvdV17}, which we pursue by employing rescaled and under-smoothing Laplace priors similar to those used in \cite{ADH21,GR22,giordano23besov,agapiou2024laplace}. Our result thus shows that such priors achieve optimal performance also in the binary classification problem. As a corollary, we obtain optimal convergence rates for the posterior mean estimator towards the spatially inhomogeneous true classification surface, see Corollary \ref{Cor:ConvPostMean}. We provide a concluding discussion in Section \ref{Sec:Discussion}, with a summary of the results and an outlook on related research questions.

%
%
%
%
%

\section{Binary classification with Besov-Laplace priors}

%
%
%

\subsection{Preliminaries and notation}\label{Sec:Notation}

Throughout, we consider functions spaces defined on the $d$-dimensional unit cube $\cube, \ d\in\N$. We write $L^p(\cube),\ p\ge 1$, for the usual Lebesgue spaces on $\cube$, equipped with norm $\|\cdot\|_2$, and $\langle\cdot,\cdot\rangle_2$ for the $L^2$-inner product.

	Let $\{\psi_{lr}, \ l\in\N, \ r=1,\dots,2^{ld}\}$ be an orthonormal basis of $L^2(\cube)$ comprising sufficiently regular, compactly supported and boundary-corrected wavelets; see e.g.~\cite[Chapter 4.3]{GN16} for details. For $\alpha\ge 0$ and $p,q\in[1,\infty]$, let the Besov space $B^\alpha_{pq}(\cube)$ be defined as
\begin{align*}
	&B^\alpha_{pq}(\cube) \\
	&\ := \left\{w \in L^p(\cube)
	 :  \|w\|_{B^\alpha_{pq}}^q := \sum_{l=1}^\infty2^{ql\left(\alpha - \frac{d}{p} + \frac{d}{2} \right)}
	\Bigg(\sum_{r=1}^{2^{ld}} |\langle w,\psi_{lr}\rangle_2|^p \Bigg)^\frac{q}{p}<\infty\right\},
\end{align*}
where the above $\ell_p$- and $\ell_q$-sequence space norms are replaced by the $\ell_\infty$-norm if $p=\infty$ or $q=\infty$ respectively. See \cite{GN16}, p.370f. Such Besov scale is known to contain the traditional Sobolev and H\"older spaces: in particular, $B^\alpha_{22}(\cube) = H^\alpha(\cube)$ for all $\alpha\ge0$, and $C^\alpha(\cube)\subseteq B^\alpha_{\infty\infty}(\cube)$ for all $\alpha>0$ (with equality holding when $\alpha\notin\N$). For $p=q=1$, the spaces $B^\alpha_{11}(\cube)$ provide a mathematical model for spatially inhomogeneous functions; in particular, the space of bounded variation functions $BV(\cube)$, which is of specific interest for many applications, enjoys the (strict) inclusions $B^1_{11}(\cube)\subset BV(\cube)\subset B^1_{1\infty}(\cube)$ (e.g.~\cite[Section 1]{DJ98}).

	When there is no risk of confusion, we omit to explicitly refer to the underlying domain, e.g.~writing $B^\alpha_{pq}$ for $B^\alpha_{pq}(\cube)$. We use the symbols $\lesssim,\ \gtrsim$, and $\simeq$, for respectively one- and two-sided inequalities holding up to universal multiplicative constants. We denote by $N(\xi;\Fcal,\delta), \ \xi>0$, for the $\xi$-covering number of a set $\Fcal$ with respect to a semi-metric $\delta$ on $\Fcal$, namely the smallest number of balls of $\delta$-radius $\xi$ required to cover $\Fcal$.

%
%
%
%
%

\subsection{Likelihood, prior and posterior}
\label{Subsec:ObsAndPrior}

Consider a random sample $(Y^{(n)},X^{(n)}):=\{(Y_i,X_i)\}_{i=1}^n$ of binary classification data arising as
\begin{equation}
\label{Eq:Obs}
\begin{split}
	Y_i|X_i &\iid \text{Bernoulli}(f(X_i)),\\
	X_i\iid \mu,
\end{split}	
\end{equation}
where $\mu$ is a given probability distribution on $\cube$ and $f:\cube\to[0,1]$ is some unknown probability response function. We are interested in the problem of estimating $f$ from data $(Y^{(n)},X^{(n)})$, whose (product) law $P^{(n)}_f$ has likelihood
\begin{equation}
\label{Eq:Likelihood}
	L_n(f):=\prod_{i=1}^n f(X_i)^{Y_i}(1-f(X_i))^{1-Y_i}.
\end{equation}
We denote by $E^{(n)}_f$ the expectation with respect to $P^{(n)}_f$. Throughout, we will make the minimal assumption that $\mu$ is absolutely continuous with respect to the Lebesgue measure and that it possesses a bounded and bounded away from zero probability density. However, we will not pursue jointly estimating $f$ and $\mu$ in the scenario where the latter is unknown.

	We are primarily interested in the setting where $f$ is spatially inhomogeneous, e.g.~possibly flat in some areas and `spiky' in others. A natural mathematical model for this scenario is then to take a (Borel measurable) parameter space $f\in\Fcal\subseteq B^\alpha_{11}(\cube)$ for some $\alpha\ge0$; see Section \ref{Sec:Notation} for definitions and details. We adopt the nonparametric Bayesian approach, endowing $f$ with the following Besov-Laplace prior. These priors represent a particular instance in the class introduced by Lassas et al.~\cite{LSS09}, and have been recently shown to achieve minimax-optimal posterior contraction rates over Besov spaces in the white noise model \cite{ADH21} as well as in other settings \cite{GR22,giordano23besov,agapiou2024laplace}. In accordance with the latter references, we construct Besov-Laplace priors for classification surfaces starting from  the random function, for $x\in [0,1]^d$,
\begin{equation}
\label{Eq:Wn}
	W_n(x) 
	= \frac{1}{n^{d/(2\alpha+d)}}
	\sum_{l=1}^\infty \sum_{r=1}^{2^{ld}} 2^{-l(\alpha-d/2)}
	W_{lr}\psi_{lr}(x), 
	\ \ W_{lr}\iid\text{Laplace},
	\ \ \alpha>d,
\end{equation}
where $\{\psi_{lr}, \ l\in\N, \ r=1,\dots,2^{ld}\}$ is an orthonormal wavelet basis of $L^2(\cube)$ as described in Section \ref{Sec:Notation}, and the Laplace (or double exponential) distribution has probability density $\lambda(z) =e^{-|z|}/2$, $z\in\R$.

	The parameter $\alpha$ in \eqref{Eq:Wn} determines the (Besov-) regularity of the realisations of the random function $W_n$: in particular, $W_n\in C(\cube)\cap B^{\alpha'}_{pq}(\cube)$ almost surely for all $\alpha'<\alpha-d$ and all $p,q\in[1,\infty]$ (cf.~Lemma 5.2 and Proposition 6.1 in \cite{ADH21}). We may then regard $W_n$ as a random element in $C(\cube)$, whose law we denote by $\Pi_{W_n}$. In accordance with the terminology of \cite{ADH21}, we call $\Pi_{W_n}$ a `rescaled $(\alpha-d)$-regular Besov-Laplace prior'.

	To construct a prior on the set of probability response  functions on $\cube$, whose range equals the compact interval $[0,1]$, we proceed transforming $W_n$ in \eqref{Eq:Wn} via the logistic link $H(z) = e^z/(e^z+1), \ z\in\R$, namely taking the random function
\begin{equation}
\label{Eq:pWn}
	f_{W_n}(x) := H[W_n(x)] = \frac{e^{W_n(x)}}{e^{W_n(x)}+1},
	\qquad x\in\cube.
\end{equation}
We denote by $\Pi_n$ the law of $f_{W_n}$ as above, which, in slight abuse of terminology, we call a rescaled $(\alpha-d)$-regular Besov-Laplace prior for the classification surface $f$.

	Given $\Pi_n$ as above, Bayes' formula (e.g.~\cite{GvdV17}, p.7) yields the posterior distribution $\Pi_n(\cdot|Y^{(n)},X^{(n)})$ of  $f|(Y^{(n)},X^{(n)})$, which is equal to
$$
	\Pi_n(A|Y^{(n)},X^{(n)}) 
	= \frac{\int_A L_n(f)d\Pi_n(f)}{\int_\Fcal L_n(f')d\Pi_n(f')},
	\qquad A\subseteq \Fcal\ \textnormal{measurable},
$$
where $L_n$ is the likelihood from \eqref{Eq:Likelihood}. According to the Bayesian paradigm, $\Pi_n(\cdot|$ $Y^{(n)},X^{(n)})$ represents the comprehensive solution to the problem of estimating $f$ from data $(Y^{(n)},X^{(n)})$ arising as in model \eqref{Eq:Obs}. It may be used to obtain point estimates (e.g.~via the posterior mean) and uncertainty quantification (via credible sets).

%
%
%
%
%

\section{Main results}
\label{Sec:MainResults}

	In this section, we study the (frequentist) consistency of $\Pi_n(\cdot|Y^{(n)},X^{(n)})$, assuming binary data $(Y^{(n)},X^{(n)})$ drawn based on some  fixed `true' unknown response probability function $f_0\in\Fcal$ to be estimated, and investigating the asymptotic concentration of the posterior distribution around $f_0$ when $n\to\infty$. We quantify the speed of concentration according to the usual notion of `posterior contraction rates', that is sequences $\varepsilon_n\to0$ such that, for large enough $M>0$,
$$
	E_{f_0}^{(n)}\Pi_n\left( f :  \delta(f,f_0)> M\varepsilon_n 
	\Big| Y^{(n)},X^{(n)}\right) \to 0
$$
as $n\to\infty$. Above, $\delta$ is a distance between classification surfaces. We refer to \cite{GvdV17} for an extensive overview on the field of the frequentist analysis of nonparametric Bayesian procedures.

	Our main result establishes that the posterior distribution resulting from rescaled $(\alpha-d)$-regular Besov-Laplace priors contracts towards ground truths $f_0\in B^\alpha_{11}(\cube)$ at minimax-optimal rate in $L^2$-distance.

\begin{theorem}\label{Theo:Main}
For fixed $\alpha>d$, let the prior $\Pi_n$ be a rescaled $(\alpha-d)$-regular Besov-Laplace prior constructed as in \eqref{Eq:pWn} for $W_n$ as in \eqref{Eq:Wn}. Let $(Y^{(n)},X^{(n)})=\{(Y_i,$ $X_i)\}_{i=1}^n$ be a random sample of binary classification data arising as in \eqref{Eq:Obs} for some fixed $f=f_0 \in B^\alpha_{11}(\cube)$ satisfying $\inf_{x\in\cube}f_0(x)>0$.  Then, for $M>0$ large enough, as $n\to\infty$,
$$
	E_{f_0}^{(n)}\Pi_n\left( f  :  \|f - f_0\|_2
	> M n^{-\frac{\alpha}{2\alpha + d}} \Big| Y^{(n)},X^{(n)}\right) \to 0.
$$
\end{theorem}

	Recall that the covariate distribution $\mu$ is assumed to posses a probability density that is bounded and bounded away from zero. In this case, the $L^2$-loss with respect to which the posterior contraction rates is obtained is equivalent to the $L^2(\mu)$-loss, and thus we could have formulated the statement of Theorem \ref{Theo:Main} in terms of the latter as well. The obtained rate $\varepsilon_n = n^{-\alpha/(2\alpha+d)}$ is known to be minimax-optimal for estimating in $L^2$-distance ground truths belonging to the space $B^\alpha_{11}(\cube)$, e.g.~\cite{DJ98}.

	Our result is in the spirit of the available theory for Gaussian priors in binary classification \cite{vdVvZ08}, for which minimax optimal posterior contraction rates in $L^2$-loss are obtained for spatially homogeneous functions under traditional Sobolev-or H\"older-type regularity assumptions. On the other hand, Gaussian priors are known to be unable to optimally model spatially inhomogeneous functions \cite{agapiou2024laplace}, and therefore they cannot be expected to achieve the same posterior contraction rates attained by Besov-Laplace priors in the setting at hand.

	The assumption that $f_0>0$ over the whole domain $[0,1]^d$ crucially guarantees that the ground truth can be written as a functional composition with respect to the employed link function $H$. Since, by construction, the support of the prior $\Pi_n$ is contained within the set of such compositions, this restriction seems unavoidable to obtain posterior consistency in the present setting. We note that the same assumption underpins the contraction rate analysis for Gaussian priors in binary classification of \cite{vdVvZ08}, where the prior is similarly constructed via compositions with (logistic-type) link functions. For the important case where we may have zero response probability at some $x\in[0,1]^d$, posterior contraction rates have been established using completely different priors constructions, in particular ones based on the Dirichlet process, e.g.~Section 9.5.6 in \cite{GvdV17}. However, we are not aware of any such result in the presence of spatially inhomogeneous classification surfaces.

	Theorem \ref{Theo:Main} establishes that the posterior distribution asymptotically concentrates over $L^2$-balls centred around the ground truth with minimal radius (in the minimax sense). The following corollary shows that the same rate of convergence $\varepsilon_n = n^{-\alpha/(2\alpha+d)}$ is attained also by the associated posterior mean, which is hence a minimax-optimal estimator of the spatially inhomogeneous classification surface. Let $\bar W_n := E^{\Pi_{W_n}}[W|Y^{(n)},X^{(n)}]$ be the mean of the posterior distribution of $W|(Y^{(n)},X^{(n)})$, arising from prior $\Pi_{W_n}$ constructed as in \eqref{Eq:pWn}. In slight abuse of terminology, let $\bar f_n:=f_{\bar W_n}=H\circ \bar W_n$, be the corresponding `posterior mean' estimate of the probability response function $f$.

\begin{corollary}\label{Cor:ConvPostMean}
Under the assumption of Theorem \ref{Theo:Main}, we have for some sufficiently large $M>0$, as $n\to\infty$,
$$
	P_{f_0}^{(n)}\left( \|\bar f_n 
	- f_0\|_2 > M n^{-\frac{\alpha}{2\alpha + d}} 
	\right)\to0.
$$
\end{corollary}

%
%
%
%
%

\section{Discussion}
\label{Sec:Discussion}

In this article, we have investigated the performance of Besov-Laplace priors for the recovery of spatially inhomogeneous binary classification surfaces. Our main result shows that, for ground truths $f_0\in B^\alpha_{11}(\cube)$, suitably calibrated rescaled $(\alpha-d)$-regular priors attain minimax-optimal posterior contraction rates in $L^2$-distance. As a corollary, we have shown that the associated posterior mean estimator is minimax-optimal. This is in line with the existing theoretical literature on Besov-Laplace priors, where similar results have been obtained in different statistical models \cite{GR22,giordano23besov,agapiou2024laplace}.

	We have focused for conciseness on non-adaptive procedures, implicitly assuming knowledge of the regularity level $\alpha$ of the ground truth for the construction of the optimally-contracting posterior distribution in Theorem \ref{Theo:Main}. Hierarchical procedures based on Besov-Laplace priors obtained by randomising the hyperparameters in the prior construction have recently been shown to achieve adaptive posterior contraction rates in the white noise model and in density estimation \cite{giordano23besov,agapiou2024adaptive}; in light of this, pursuing similar approaches for the binary classification problem is likely feasible.

	In Theorem \ref{Theo:Main}, minimax-optimal posterior contraction rates for ground truths $f_0\in B^\alpha_{11}(\cube)$ have been obtained by rescaled $(\alpha-d)$-regular Besov-Laplace priors. The same specific prior tuning, that finely exploits the interplay between rescaling and undersmoothing, has been employed in the other available results for Besov-Laplace priors \cite{GR22,giordano23besov,agapiou2024laplace}, and marks an interesting difference with the existing literature on Gaussian priors, where minimax-optimal posterior contraction rates are typically obtained by smoothness-matching (and non-rescaled) priors, in the presence of spatially homogeneous ground truths. This is issue is extensively discussed in \cite{giordano23besov}. Recently, Dolera et al.~\cite{dolera2024strong} have shown, in the white noise model, that non-rescaled smoothness-matching Besov-Laplace priors optimally contract towards ground truths belonging to the Besov scale $B^\alpha_{11}$. However, it remains unclear at present whether their result (and proof techniques) can be extended to the setting at hand.

	Lastly, the computational aspects of posterior inference with Besov-Laplace priors in binary classification represent a further interesting avenue for future research. These priors naturally lend themselves to discretisation by truncating the series in \eqref{Eq:Wn} at some user-specified (sufficiently high) level $L\in\N$, giving rise to high-dimensional log-concave product distributions. Combined with the simple likelihood structure from \eqref{Eq:Likelihood}, this furnishes a basis for the implementation of Markov chain Monte Carlo (MCMC) methods to sample from the resulting posterior distribution. For recent developments concerning MCMC algorithms suited Besov-Laplace priors in infinite-dimensional statistical models, we refer to \cite{chen2018dimension}, where more references can be found.
	
%
%
%
%
%

\section{Proofs}
\label{Sec:Proof}

%
%
%
%
%

\subsection{Proof of Theorem \ref{Theo:Main}}

Since $f_0$ is strictly positive, we may write $f_0 = f_{w_0}:=H\circ w_0$ for $H$ the logistic link function and $w_0:= H^{-1}\circ f_0$. In particular, as $f_0\in B^\alpha_{11}$ and $H^{-1}$ is smooth, we have by Theorem 10 in \cite{BS10} that $w_0\in B^\alpha_{11}$.

	The proof of Theorem \ref{Theo:Main} is based on the general theory for posterior contraction rates in i.i.d.~sampling models, adapted to the case of binary classification data and nonparametric priors constructed via the logistic link function. In particular, combining Theorem 8.9 in \cite{GvdV17}  (and its proof) with Lemma 2.8 in the same reference shows that if for some sequence $\varepsilon_n\to0$ such that $n\varepsilon_n^2\to\infty$, some constant $C>0$, and all $n\in\N$ large enough,
\begin{equation}
\label{Eq:SmallBall}
	\Pi_{W_n}\left(w  :  \|w - w_0\|_{L^2(\mu)}\le \varepsilon_n\right)
	\ge  e^{-Cn\varepsilon_n^2},
\end{equation}
and if there exists measurable sets $\Wcal_n\subseteq C(\cube)$ such that
\begin{equation}
\label{Eq:Sieves}
	\Pi_{W_n}(\Wcal_n^c)\le e^{-(C+4)n\varepsilon_n^2};
	\qquad
	\log N(\varepsilon_n; \Wcal_n, \|\cdot\|_{L^2(\mu)})
	\lesssim n\varepsilon_n^2,
\end{equation}
then, for sufficiently large $M>0$ and some $\bar C>0$, as $n\to\infty$,
\begin{equation}
\label{Eq:GeneralConclusion}
	E_{f_0}^{(n)}
	\Pi_n\left(f_w: w\in\Wcal_n,\ \|f_w - f_0\|_{L^2(\mu)}\le M\varepsilon_n\big|
	Y^{(n)},X^{(n)}\right)
	\ge 1 - e^{-\bar C n\varepsilon_n^2}.
\end{equation}

	We proceed verifying the `small ball probability lower bound' \eqref{Eq:SmallBall} and the `sieve condition' \eqref{Eq:Sieves}, wherein we note that we may replace the norm $\|\cdot\|_{L^2(\mu)}$ with $\|\cdot\|_{2}$ since $\mu$ has a bounded and bounded away from zero probability density. Set $\varepsilon_n := n^{-\alpha/(2\alpha + d)}$, and write $W_n$ in \eqref{Eq:Wn} as $W_n = (n\varepsilon_n^2)^{-1}W$, where $W$ is the non-rescaled $(\alpha-d)$-regular Besov-Laplace random element
\begin{equation}
\label{Eq:W}
	W(x) := 
	\sum_{l=1}^\infty \sum_{r=1}^{2^{ld}} 2^{-l\left(\alpha - d/2\right)}
	W_{lr}\psi_{lr}(x),
	\qquad x\in\cube,
	\qquad W_{lr}\iid \textnormal{Laplace}.
\end{equation}
Writing $\Pi_{W}$ for the law of $W$, we then have
$$
	\Pi_{W_n}( w  :  \|w - w_0\|_2 \le  \varepsilon_n )
	=\Pi_{W}\left(w: \|w - n\varepsilon_n^2w_0\|_2 \le  n \varepsilon_n^3
	\right )
$$
so that by using the `decentering inequality' and the sup-norm small ball probability lower bound from eq.~(32) and (34) in \cite{giordano23besov} respectively, we obtain that the probability of interest is greater than
\begin{align*}
	e^{-\|w_0\|_{B^\alpha_{11}}n\varepsilon_n^2}
	\Pi_{W}\left(w:\|w\|_\infty\le c_1 n \varepsilon_n^3\right)
	&\ge 
	e^{-\|w_0\|_{B^\alpha_{11}}n\varepsilon_n^2}
	\Pi_{W}\left(w:\|w\|_\infty\le c_1 n \varepsilon_n^3\right)\\
	&
		\ge e^{-\|w_0\|_{B^\alpha_{11}}n\varepsilon_n^2+c_2(n\varepsilon_n^3)^{-d/(\alpha - d)}}\\
	&=
		e^{-Cn\varepsilon_n^2},
\end{align*}
for constants $c_1,c_2,C>0$. This concludes the verification of \eqref{Eq:SmallBall}. Turning to condition \eqref{Eq:Sieves}, define the sieves
\begin{align}
\label{Eq:Wcaln}
	\Wcal_n :=\left\{w = w^{(1)} + w^{(2)} : \|w^{(1)}\|_2\le R \varepsilon_n, 
	\|w^{(2)}\|_{B^\alpha_{11}}\le R  ,\ \|w\|_\infty\le R\right\}.
\end{align}
Then
\begin{align*}
	\Pi_{W_n}(\Wcal_n^c)
	&\le 2-\Pi_{W_n}\Big( w= w^{(1)} + w^{(2)} 
	: \|w^{(1)}\|_2\le R \varepsilon_n,  \ 
	\|w^{(2)}\|_{B^\alpha_{11}}\le R  \Big)\\
	&\quad\  - \Pi_{W_n}(w: \|w\|_\infty\le R ),
\end{align*}
which, by Lemma 9 in \cite{giordano23besov}, fixing $K>C+4$ and choosing sufficiently large $R>0$, is smaller than
$$
	2e^{-Kn\varepsilon_n^2}\le e^{-(C+4)n\varepsilon_n^2}.
$$
This proves the first inequality in \eqref{Eq:Sieves}. For the second, note that since $\Wcal_n \subset \left\{ w = w^{(1)} + w^{(2)} : \|w^{(1)}\|_2\le R \varepsilon_n,\  
\|w^{(2)}\|_{B^\alpha_{11}}\le R  \right\}$, using the embedding of $B^\alpha_{11}$ into $B^\alpha_{1\infty}$ (e.g., Theorem 3.3.1 in \cite{T83}) and the complexity bound from Theorem 4.3.36 in \cite{GN16}, we have
$$	
	\log N(\varepsilon_n; \Wcal_n, \|\cdot\|_2)
	\lesssim
	\log N(\varepsilon_n; \{w :  \|w\|_{B^\alpha_{1\infty}}\le R  \}, \|\cdot\|_2)
	\lesssim \varepsilon_n^{-\frac{d}{\alpha}} = n\varepsilon_n^2,
$$
concluding the verification of \eqref{Eq:Sieves} and, in view of \eqref{Eq:GeneralConclusion}, the proof of Theorem \ref{Theo:Main}.\qed

%
%
%
%
%

\subsection{Proof of Corollary \ref{Cor:ConvPostMean}}

Inspection of the proof of Theorem \ref{Theo:Main}, cf.~\eqref{Eq:GeneralConclusion}, shows that it holds that for some sufficiently large $M>0$,
$$
	E_{f_0}^{(n)}
	\Pi_n\left(f_w: w\in\Wcal_n,\ \|f_w - f_0\|_2\le M\varepsilon_n\big|
	Y^{(n)},X^{(n)}\right)
	\ge 1 - e^{-\bar C n\varepsilon_n^2},
$$
where $\Wcal_n$ is the sieve defined in \eqref{Eq:Wcaln}. Since the latter is included in a sup-norm ball and the logistic link function $H$ is continuous with continuous inverse $H^{-1}$, it then follows that
\begin{equation}
\label{Eq:StrongThesis}
	E_{f_0}^{(n)}
	\Pi_{W_n}\left(w\in\Wcal_n: \|w - w_0\|_2\le M \varepsilon_n\big|
	Y^{(n)},X^{(n)}\right)
	\ge 1 - e^{-\bar C n\varepsilon_n^2},
\end{equation}
having denoted by $\Pi_{W_n}(\cdot|Y^{(n)},X^{(n)})$ the posterior distribution of $W|(Y^{(n)},X^{(n)})$ arising from the prior $\Pi_{W_n}$ constructed as in \eqref{Eq:Wn}. Write shorthand $D^{(n)}:=(Y^{(n)},X^{(n)})$. For $\bar W_n = E^{\Pi_{W_n}}[W|D^{(n)}]$ the posterior mean, we proceed showing that, as $n\to\infty$,
$$
	P_{f_0}^{(n)}\left( \|\bar W_n 
	- w_0\|_2 > M n^{-\frac{\alpha}{2\alpha + d}} 
	\right)\to0.
$$
Upon doing so, the conclusion of Corollary \ref{Cor:ConvPostMean} directly follows since $H$ is Lipschitz. We use uniform integrability arguments for Besov-Laplace priors developed in \cite{GR22}. Start by noting that by Jensen's inequality
\begin{align*}
	\| \bar W_n - w_0\|_2
	\le
		E^{\Pi_{W_n}}\left[ \| W - w_0\|_2 \big|D^{(n)}\right],
\end{align*}
which, for any $M>0$, using the Chauchy-Schwarz inequality, is upper bounded by
\begin{align*}
	&
		M\varepsilon_n +
		E^{\Pi_{W_n}}\left[ \| W - w_0\|_2 1_{\{\| W - w_0\|_2 >M\varepsilon_n\}}
		\big|D^{(n)}\right] \\
	&\quad \le		
			M\varepsilon_n + 	E^{\Pi_{W_n}}\left[ \| W - w_0\|_2 ^2\big|D^{(n)}\right]^{1/2}
	\Pi_{W_n}(w : \| w - w_0\|_2 >M\varepsilon_n|D^{(n)})^{1/2}.
\end{align*}
Using the above and \eqref{Eq:StrongThesis}, for any $A>0$, we can then find $M>0$ large enough such that
\begin{align*}
	& P_{f_0}^{(n)} (\| \bar W_n - w_0\|_2> 2M\varepsilon_n) \\
	& \leq P_{f_0}^{(n)} \left( E^{\Pi_{W_n}}\left[ \| W - w_0\|_2 ^2\big|D^{(n)}\right]^\frac{1}{2}
		\Pi_{W_n}(w : \| w - w_0\|_2 >M\varepsilon_n|D^{(n)})^\frac{1}{2} > M\varepsilon_n\right)\\
	& \leq P_{f_0}^{(n)} \left( E^{\Pi_{W_n}}\left[ \| W - w_0\|_2 ^2\big|D^{(n)}\right]
		e^{-An\varepsilon_n^2} > M^2 \varepsilon_n^2 \right) + o(1).
\end{align*}
For $L_n$ the likelihood from \eqref{Eq:Likelihood} and $C>0$ the constant in the small ball estimate \eqref{Eq:SmallBall}, define the event
$$
	A_n = \left\{ \int_{C(\cube)}
	L_n(f_w)/L_n(f_0)d\Pi_{W_n}(w) \geq e^{-(C+2)n\varepsilon_n^2} \right\},
$$
which, by Lemma 8.10 in \cite{GvdV17}, in view of the validity of \eqref{Eq:SmallBall}, satisfies $P_{f_0}^{(n)}(A_n)\to 1$ as $T\to\infty$. The second probability in the second to last display is then upper bounded by
\begin{align*}
	&P_{f_0}^{(n)}\left( \left\{ \frac{\int_{C(\cube)}\| w - w_0\|^2_2 
	L_n(f_w)/L_n(f_0)d\Pi_{W_n}(w)} 
	{\int_{C(\cube)}
	L_n(f_{w'})/L_n(f_0)d\Pi_{W_n}(w')}
	>M^2 \varepsilon_n^2 e^{An\varepsilon_n^2} \right\} \cap  A_n\right)\\
	& \quad \le
	P_{f_0}^{(n)}\left(\int_{C(\cube)}\| w - w_0\|^2_2 
	\frac{L_n(f_w)}{L_n(f_0)}d\Pi_{W_n}(w)
	>M^2 \varepsilon_n^2 e^{(A-C-2)n\varepsilon_n^2}\right) + o(1).
\end{align*}
By Markov's inequality and Fubini's theorem, the latter probability is smaller than
\begin{align*}
		\frac{e^{-(A-C-2)n\varepsilon_n^2}}{M^2 \varepsilon_n^2} 
		 \int_{C(\cube)}\| w - w_0\|^2_2 &
		E_{f_0}^{(n)}  \left[ \frac{L_n(f_w)}{L_n(f_0)} \right] d\Pi_{W_n}(w) \\
		&\qquad =
		  \frac{e^{-(A-C-2)n\varepsilon_n^2}}{M^2 \varepsilon_n^2} E^{\Pi_{W_n}}\left[ \| W - w_0\|^2_2\right] ,
\end{align*}
since by a change of measure $E_{f_0}^{(n)}\left[L_n(f_w)/L_n(f_0)\right]=E_{f_w}^{(n)}[1]=1$. To proceed, we note that the Fernique-like theorem in Section 2 of \cite{DHS12} implies that
\begin{align*}
	E^{\Pi_{W_n}} [\|W\|_2 ^2]
	\lesssim 1.
\end{align*}
Combining the last three displays then yields that for any $A>0$, we can find $M>0$ large enough such that  
\begin{align*}
	P_{f_0}^{(n)} (\| \bar W_n - w_0\|_2> 2M\varepsilon_n) 
	\lesssim \frac{e^{-(A-C-2)n\varepsilon_n^2}}
	{M^2 \varepsilon_n^2} + o(1)
\end{align*}
since $w_0\in B^\alpha_{11}\subset L^2$ as $\alpha>d$. Taking $A>(C+2)$ and using that $\varepsilon_n = n^{-\alpha/(2\alpha+d)}$, the right-hand side tends to zero as $n \to\infty$.\qed

%
%
%
%
%

\paragraph{Acknowledgements.}

We are grateful to one anonymous referee for helpful comments that lead to an improvement of this article. This research was supported by the Italian Ministry of Education, University and Research (MIUR), “Dipartimenti di Eccellenza” grant 2023-2027, and partially by MIUR, PRIN project 2022CLTYP4.

\bibliographystyle{acm}

\bibliography{References}

\end{document}